\documentclass[a4paper]%
              {amsart}

\renewcommand\newsymbol[5]{%
\DeclareMathSymbol#1{#3}{\ifcase #2\or AMSa\or AMSb\fi}{"#4#5}}

\newcommand\newbbbletter[2]{%
\DeclareMathSymbol#1{0}{AMSb}{`#2}}

\newbbbletter\C    C
\newbbbletter\Kat  K    
\newbbbletter\LL   L
\newbbbletter\M    M    
\newbbbletter\N    N    
\newbbbletter\Po   P    
\newcommand\bound{\mathfrak{b}}
\newcommand\cont{\mathfrak{c}}
\newsymbol\forces  130D
\newsymbol\restr   1216
\newsymbol\le      1236
\newsymbol\ge      123E
\newsymbol\symmdif 124D
\newsymbol\sulk    1061
\newsymbol\nsubseteq 232A
\let\diamond\diamondsuit
\DeclareMathSymbol\adA{0}{symbols}{`A}  
\DeclareMathSymbol\calB{0}{symbols}{`B}
\DeclareMathSymbol\calC{0}{symbols}{`C}
\DeclareMathSymbol\calF{0}{symbols}{`F}
\DeclareMathSymbol\calG{0}{symbols}{`G}
\DeclareMathSymbol\idI{0}{symbols}{`I}  
\DeclareMathSymbol\calL{0}{symbols}{`L}
\DeclareMathSymbol\Null{0}{symbols}{`N} 
\DeclareMathSymbol\Meag{0}{symbols}{`M} 
\DeclareMathSymbol\pow{0}{symbols}{`P}  
\let\savephi\phi
\let\phi\varphi
\let\varphi\savephi
\newcommand\Bor{\operatorname{Bor}}

\newcommand\Fn{\operatorname{Fn}}
\newcommand\concat{\mathord{{^\sulk}}}
\newcommand\fin{\mathit{fin}}
\newcommand\functions{{}^\omega\omega}
\newcommand\bintree{2^{<\omega}}
\newcommand\sbintree[1][f]{2^{#1(n)}}
\newcommand\almeq{\mathrel{=^*}}
\newcommand\almneq{\mathrel{\neq^*}}
\newcommand\almsubseteq{\mathrel{\subseteq^*}}
\newcommand\alml{\mathrel{<^*}}
\newcommand\almle{\mathrel{\le^*}}
\newcommand\inv{^\leftarrow}
\newcommand\powN{\pow(\N)}
\newcommand\powNfin{\powN/\fin}
\newcommand\betaNminN{\beta\N\setminus\N}
\newcommand\orpr[2]{\langle#1,#2\rangle}

\newcommand\card[1]{\lvert#1\rvert}
\newcommand\bigcard[1]{\bigl|#1\bigr|}
\newcommand\sPhin[1][n]{\Phi\bigl(\{#1\}\bigr)}
\newcommand\axiom[1]{\mathsf{#1}}
\newcommand\CH{\axiom{CH}}
\newcommand\MA{\axiom{MA}}
\newcommand\OCA{\axiom{OCA}}
\newcommand\PFA{\axiom{PFA}}
\newcommand\ZFC{\axiom{ZFC}}
\newcommand\Cantor{2^\omega}
\newcommand\szer{s\concat0}
\newcommand\sone{s\concat1}
\newcommand\Cohen{\Fn(\omega,2)}
\def\iff/{if\kern0ptf}
\newtheorem*{maintheorem}{Main Theorem}
\newtheorem{theorem}{Theorem}[section]
\newtheorem{proposition}[theorem]{Proposition}
\newtheorem{lemma}[theorem]{Lemma}
\theoremstyle{definition}
\newtheorem{question}{Question}[section]

\newcounter{claim}[theorem]
\renewcommand\theclaim{\arabic{claim}}
           {\endgraf\stepcounter{claim}\smallskip
            
            \noindent\textsc{Claim \theclaim}. \ignorespaces}%
           {\parfillskip0pt\hfil$\symmdif$\endgraf\smallskip}

\newcommand\cprime{$'$}

\begin{document}

\title{The measure algebra does not always embed}

\author[Alan Dow]{Alan Dow*}
\thanks{*The research of the first author was supported in part by
        the Netherlands Organization for Scientific Research (NWO) ---
        Grant~B\,61-408.}
\address{Department of Mathematics\\
         York University\\
         4700 Keele Street\\
         Toronto, Ontario\\
         Canada M3J 1P3}
\email{dowa@YorkU.CA}
\author{Klaas Pieter Hart}
\address{Faculty of Information Technology and Systems\\
         TU Delft\\
         Postbus 5031\\
         2600~GA {} Delft\\
         the Netherlands}
\email{k.p.hart@twi.tudelft.nl}
\urladdr{http://aw.twi.tudelft.nl/\~{}hart}

\subjclass{Primary: 28A60. Secondary: 06E99, 03E35, 54G05}
\keywords{measure algebra, embedding, Open Colouring Axiom, $\powNfin$}

\begin{abstract}
The Open Colouring Axiom implies that the measure algebra cannot be
embedded into~$\powNfin$.
We also discuss errors in previous results on the embeddability
of the measure algebra.
\end{abstract}

\maketitle

\section*{Introduction}

The aim of this paper is to prove the following result.

\begin{maintheorem}
The Open Colouring Axiom implies that the measure algebra cannot be
embedded into the Boolean algebra~$\powNfin$.
\end{maintheorem}

By `the measure algebra' we mean the quotient of the $\sigma$-algebra
of Borel sets of the real line by the ideal of sets of measure zero.

There are various reasons, besides sheer curiosity, why it is of interest
to know whether the measure algebra can be embedded into~$\powNfin$.
One reason is that there is great interest in determining what the
subalgebras of~$\powNfin$ are.
One of the earliest and most influential result in this direction is
Parovi\v{c}enko's theorem from~\cite{Parovicenko63}, which states that every
Boolean algebra of size~$\aleph_1$ can be embedded into~$\powNfin$ ---
with the obvious corollary that the Continuum Hypothesis ($\CH$) implies
that $\powNfin$ is a universal Boolean algebra of size~$\cont$:
a Boolean algebra embeds into $\powNfin$ \iff/ it is of size~$\cont$ or less.
It is therefore natural to ask how much of this universality remains
without assumptions beyond~$\ZFC$.
It has long been known that every $\sigma$-centered Boolean algebra embeds
into~$\powNfin$ but the question for more general c.c.c.\ Boolean
algebras has proven to be much more difficult ---
with the case of the measure algebra being seen as a touchstone.

This particular case was especially interesting since, by Stone duality,
an embedding of the measure into~$\powNfin$ would provide a
c.c.c.\ nonseparable remainder of~$\N$.
A $\ZFC$ construction of such a remainder was given by~Bell
in~\cite{Bell80}; such remainders were put to good use by Van~Mill
in~\cite{vanMill82}, see also his survey~\cite{vanMill84}.
The question of the embeddability of the measure algebra remained however.

In recent years many results about maps between Boolean algebras or
topological spaces, which were shown to hold under~$\CH$, were shown
to fail under~$\OCA$ --- see for example
\cite{Just92b,Farah97,DowHart99}.
Our result and its proof fall into the same category: in all cases
$\OCA$~implies that the desired map should have some simple structure
whereas one can show, usually in~$\ZFC$, that the desired map cannot have
this simple structure.

We should note that our result clashes with the result of Frankiewicz
and Gutek from~\cite{FrankiewiczGutek81}, which says that Martin's Axiom
implies that the measure algebra can be embedded into~$\powNfin$;
for completeness we shall point out a gap in their proof.
In addition we feel that we should mention the paper~\cite{Frankiewicz85}.
The main result in that paper is that one can establish the consistency of
the nonembeddability of the measure algebra using Shelah's oracle-c.c.\
method.
Regrettably, again the argumentation appears to be incomplete as we shall
discuss further in Section~\ref{sec.comments}.

The paper is organized as follows.
Section~\ref{sec.prelim} contains the necessary preliminaries, including
a discussion of what `simple structure' means in our context.
Section~\ref{sec.exact.is.not.complete} shows that no embedding
can have this simple structure.
In Sections~\ref{sec.OCA.gives.exact.and.complete}
and~\ref{sec.powN.into.powNfin} we show how $\OCA$ implies that
embeddings must have a simple structure.
Finally then Section~\ref{sec.comments} is devoted to the discussion
of problems in the previously published work in this area. This will
allow us to resurrect the interesting question of whether this result
can be established with the oracle-c.c.\ method.

\section{Preliminaries}
\label{sec.prelim}

\subsection*{The Measure Algebra}
The standard representation of the Measure Algebra is as the quotient of
the $\sigma$-algebra of Borel sets of the unit interval by the ideal
of sets of Lebesgue measure zero.
For ease of notation we choose a different underlying set, namely
$\C=\omega\times\Cantor$, where $\Cantor$is the Cantor set.
We consider the Cantor set endowed with the natural coin-tossing
measure~$\mu$, determined by specifying $\mu\bigl([s]\bigr)=2^{-\card{s}}$.
Here $s$ denotes a finite partial function from~$\omega$ to~$2$
and $[s]=\{x\in\Cantor:s\subset x\}$.
We extend $\mu$ on the Borel sets of~$\C$ by setting
$\mu\bigl(\{n\}\times[s]\bigr)=2^{-\card{s}}$ for all~$n$ and~$s$.

The measure algebra is isomorphic to the quotient algebra
$\M=\Bor(\C)/\Null$, where $\Null=\{N\subseteq\C:\mu(N)=0\}$;
henceforth we shall work with $\M$.

\subsection*{Liftings of embeddings}

Assume $\phi:\M\to\powNfin$ is an embedding of Boolean algebras and take a
\emph{lifting} $\Phi:\M\to\powN$ of~$\phi$; this is a map that chooses a
representative~$\Phi(a)$ of~$\phi(a)$ for every~$a$ in~$\M$.

We shall be working mostly with the restrictions of~$\phi$ and~$\Phi$ to the
family of (equivalence classes of) open subsets of~$\C$ and in particular with
their restrictions to the canonical base for~$\C$, which is
$$
\calB=\bigl\{\{n\}\times[s]:n\in\omega, s\in\bintree\bigr\}.
$$
To keep our formulas manageable we shall identify $\calB$ with the
set $B=\omega\times\bintree$.
We shall also be using layers/strata of~$B$ along functions from~$\omega$
to~$\omega$: for $f\in\functions$ we put
$B_f=\bigl\{\orpr{n}{s}:n\in\omega, s\in\sbintree\bigr\}$.

For a subset $O$ of~$B_f$ we abbreviate
$\phi\bigl(\bigcup\bigl\{\{n\}\times[s]:\orpr{n}{s}\in O\bigr\}\bigr)$
by~$\phi(O)$ and define $\Phi(O)$ similarly.
Observe that $O\mapsto\phi(O)$ defines an embedding of~$\pow(B_f)$
into~$\powNfin$.
As an extra piece of notation we use $\Phi[O]$ (square brackets) to denote the
union $\bigcup\bigl\{\Phi(n,s):\orpr{n}{s}\in O\bigr\}$,
where $\Phi(n,s)$ abbreviates $\Phi\bigl(\bigl\{\orpr{n}{s}\bigr\}\bigl)$.

For later use we explicitly record the following easy lemma.

\begin{lemma}\label{lemma.Phi.fin.add}
If $f\in\functions$ and if\/ $O$ is a finite subset of\/~$B_f$ then
$\Phi(O)\almeq\Phi[O]$.
\end{lemma}

\begin{proof}
Both sets represent $\phi(O)$.
\end{proof}

Let us call a lifting \emph{complete} if it satisfies
Lemma~\ref{lemma.Phi.fin.add} for \emph{every}~$f\in\functions$ and
\emph{every} subset~~$O$ of~$B_f$.

We can always make a lifting~$\Phi$ \emph{exact}, by which we mean that
the sets~$\Phi(n,\emptyset)$ form a partition of~$\N$ and that
every~$\Phi(s,n)$ is the disjoint union of~$\Phi(n,\szer)$
and~$\Phi(n,\sone)$;
indeed, we need only change each of the countably many sets~$\Phi(n,s)$ by
adding or deleting finitely many points to achieve this.

Our proof may now be summarized in a few lines:

\begin{enumerate}
\item For every exact lifting~$\Phi$ of an embedding~$\phi$ there are
      an~$f\in\functions$ and an infinite subset~$O$ of~$B_f$ such that
      $\Phi(O)\almneq\Phi[O]$, i.e., no exact lifting is complete
      --- see Proposition~\ref{prop.exact.is.not.complete}.
\item $\OCA$~implies that every embedding~$\phi$ gives rise to an
      embedding~$\varphi$ with a lifting~$\Phi$ that is both exact and
      complete --- see Section~\ref{sec.OCA.gives.exact.and.complete}.
\end{enumerate}

Together these two statements show that under $\OCA$ there
cannot be any embedding of $\M$ into~$\powNfin$ at all.

\subsection*{The Open Colouring Axiom}

The \emph{Open Coloring Axiom} ($\OCA$) was formulated by
To\-dor\-\v{c}e\-vi\'c in~\cite{Todorcevic89}.
It reads as follows:
if $X$ is separable and metrizable and if $[X]^2=K_0\cup K_1$,
where $K_0$~is open in the product topology of~$[X]^2$, then
\emph{either} $X$~has an uncountable $K_0$-homogeneous subset~$Y$
\emph{or} $X$~is the union of countably many $K_1$-homogeneous subsets.

One can deduce $\OCA$ from the \emph{Proper Forcing Axiom} ($\PFA$)
or prove it consistent in an $\omega_2$-length countable support
proper iterated forcing construction, using~$\diamond$ on~$\omega_2$
to predict all possible subsets of the Hilbert cube and all possible
open colourings of these.

We shall not apply $\OCA$ directly but use some of its known consequences
to prove our main result.
A major application occurs in Section~\ref{sec.powN.into.powNfin},
where we rely on a result from~\cite{DowHart99} regarding the behaviour
of embeddings of~$\powNfin$ into itself.

Here and in the next subsection we collect some results of a more general
nature.
To begin a definition:
$\bound$~is the minimum cardinality of a family~$\calF$ of functions
from~$\N$ to~$\N$ for which there is no upper bound with respect
to~$\almle$,
i.e., whenever $g\in\N^\N$ there if~$f\in\calF$ such that
$\bigl\{n:f(n)>g(n)\bigr\}$ is infinite.

The first consequence of~$\OCA$ that we need is the equality
$\bound=\aleph_2$; it was established in~\cite[Theorem~3.16]{Bekkali91}.
It follows that the following lemma also holds under~$\OCA$.

\begin{lemma}[$\bound\ge\aleph_2$]\label{lemma.eventually.constant}
Assume $f\mapsto\alpha_f$ is a map from~$\functions$ to~$\omega_1$
that is monotone with respect to~$\almle$ and $\in$, i.e.,
if $f\almle g$ then $\alpha_f\le\alpha_g$.
Then the map is bounded, i.e., there is an~$\alpha$ such
that $\alpha_f\le\alpha$ for all~$f$.
\end{lemma}

\begin{proof}
Because $\bound\ge\aleph_2$ there is a~$\almle$-cofinal family~$\calF$
in~$\functions$ on which our map constant, say with value~$\alpha$.
Because the map is monotone this~$\alpha$ is the ordinal that we are looking
for.
\end{proof}

\subsection*{Coherent families of functions}

Twice in our proof we shall want to combine a family of partial functions
into one single function.
In both cases we shall have an ideal~$\idI$ of subsets of some countable
set~$C$ and for each~$I$ a function~$f_I$ with domain~$I$ such that
whenever $I\subseteq J$ in~$\idI$ we have $f_J\restr I\almeq f_I$
--- such a family will be called \emph{coherent}.
The following theorem, which is Theorem~3.13 from~\cite{Bekkali91},
tells us when a coherent family can be \emph{uniformized},
i.e., when we can get one function~$f$ with domain~$\bigcup\idI$ such that
$f\restr I\almeq f_I$ for all~$I$.

\begin{theorem}[$\OCA$]\label{thm.uniformizing}
If $\idI$ is a $P_{\aleph_1}$-ideal then every coherent family of functions
on~$\idI$ with values in~$\omega$ can be uniformized.
\end{theorem}

An ideal $\idI$ is a $P_{\aleph_1}$-ideal if for every subfamily~$\idI'$
of~$\idI$ of size~$\aleph_1$ (or less) one can find an element~$J$ of~$\idI$
such that $I\almsubseteq J$ for all~$I\in\idI'$.

We shall need the following generalization of Theorem~\ref{thm.uniformizing}
--- it actually turns out to be a special case.

\begin{theorem}[$\OCA$]\label{cor.uniformizing}
If $\idI$ is a $P_{\aleph_1}$-ideal on~$\omega$ then every coherent family of
functions on~$\idI$ with values in~$\powN$ can be uniformized.
\end{theorem}

\begin{proof}
Let $\{f_I:I\in\idI\}$ be a coherent family of functions, with values
in~$\powN$.
For $I\in\idI$ and $g\in\functions$ let
$L_{I,g}=\bigl\{\orpr{n}{m}:n\in I, m\le g(n)\bigr\}$
and $R_{I,f}=\bigl\{\orpr{n}{m}:n\in I,m\le g(n), m\in f_I(n)\bigr\}$.
The sets $L_{I,g}$ generate a $P_{\aleph_1}$-ideal on the
countable set~$\N\times\N$ and one readily checks that
$R_{I,g}\almeq R_{J,h}\cap L_{I,g}$ whenever $I\subseteq J$ and $g\alml h$.

Now apply Theorem~\ref{thm.uniformizing} to find $R\subseteq\N\times\N$
such that $R\cap L_{I,g}\almeq R_{I,g}$ for all~$I$ and~$g$.
This defines $f:\omega\to\powN$ by $m\in f(n)$ iff $\orpr{n}{m}\in R$.

If $f\restr I$ were not almost equal to~$f_I$ then we'd find infinitely
many~$n$ with an~$m_n$ in~$f(n)\symmdif f_I(n)$.
But then $R\cap L_{I,g}\almneq R_{I,g}$, where $g\in\functions$ follows
$n\mapsto m_n$.
\end{proof}

\section{No exact lifting is complete}
\label{sec.exact.is.not.complete}

Assume $\phi:\M\to\powNfin$ is an embedding and consider an exact lifting~$\Phi$
of~$\phi$.
The following proposition shows that $\Phi$~is not complete.

\begin{proposition}\label{prop.exact.is.not.complete}
There is a sequence $\langle t_n:n\in\omega\rangle$ in~$\bintree$
such that for the open set $O=\bigcup_{n\in\omega}\{n\}\times[t_n]$ we have
$\Phi(O)\almneq\Phi[O]$.
\end{proposition}

\begin{proof}
Take, for each~$n$, the monotone enumeration~$\{k(n,i):i\in\omega\}$
of~$\Phi(n,\emptyset)$ and apply the equalities to find $t(n,i)\in2^{i+2}$
such that $k(n,i)\in\Phi\bigl(n,t(n,i)\bigr)$.
Use these~$t(n,i)$ to define open sets
$U_n=\bigcup_{i\in\omega}\{n\}\times\bigl[t(n,i)\bigr]$;
observe that $\mu(U_n)\le\sum_{i\in\omega}2^{-i-2}=\frac12$.
It follows that $\Phi\bigl(\{n\}\times U_n^c\bigr)$ is infinite.

We let $F$ be the closed set $\bigcup_{n\in\omega}\{n\}\times U_n^c$;
its image~$\Phi(F)$ meets every~$\Phi(n,\emptyset)$ in an infinite set.
For every~$n$ let $i_n$ be the first index with $k(n,i_n)\in\Phi(F)$
and consider the open set
$O=\bigcup_{n\in\omega}\{n\}\times\bigl[t(n,i_n)\bigr]$
and the infinite set~$I=\bigl\{k(n,i_n):n\in\omega\bigr\}$.

Observe the following
\begin{enumerate}
\item $\Phi(O)\cap\Phi(F)\almeq\emptyset$, because $O\cap F=\emptyset$;
\item $I\subseteq\Phi(F)$, by our choice of the $i_n$; and
\item $I\subseteq\Phi[O]$, by the choice of the $t(n,i_n)$.
\end{enumerate}
It follows that $\bigl<t(n,i_n):n\in\omega\bigr>$ is as required.
\end{proof}

This completes the proof of the first half of the main argument.

\section{$\OCA$ gives embeddings with exact liftings that are complete}
\label{sec.OCA.gives.exact.and.complete}

We assume $\phi$~is an embedding of $\M$ into~$\powNfin$.
We shall find, assuming~$\OCA$, an infinite set~$A$ and a lifting~$\Phi$
of~$\phi$ that is exact and complete when restricted to subsets
of~$A\times\bintree$.

The infinite set~$A$ will come from an almost disjoint family
on~$\omega$: we fix bijection~$\pi$ between~$\omega$ and~$\bintree$ and define
for~$x\in\Cantor$ the set~$A_x$
by $A_x=\pi\inv\bigl[\{x\restr n:n\in\omega\}\bigr]$.

For the rest of this section we fix an $\aleph_1$-sized subfamily~$\adA$
of the~$A_x$'s and enumerate it as $\{A_\alpha:\alpha<\omega_1\}$.
Using~$\OCA$ we shall show that all but countably many~$A_\alpha$ are as
required.

By construction the family~$\adA$ has a special property, commonly referred
to as \emph{neatness}; an almost disjoint family~$\calC$ is neat if there
is a bijection~$\pi$ between~$\omega$ and~$\bintree$ such that for
every~$C\in\calC$ the set~$\pi[C]$ is a subset of some branch~$x_C$ and,
moreover, the map $C\mapsto x_C$ is one-to-one.

Our final lifting will be a limit, via Theorem~\ref{cor.uniformizing}, of a
coherent family of liftings; these liftings will be defined only partially
so we fix an exact lifting~$\Psi$ of~$\phi$ to extend these partial liftings.

The key technical result is the following; we postpone its proof
until Section~\ref{sec.powN.into.powNfin}.

\begin{theorem}[$\OCA$]\label{thm.powN.into.powNfin}
Let $\phi$ be an embedding of~$\powN$ into $\powNfin$ and let
$\adA$ be a neat almost disjoint family on~$\N$ of size~$\aleph_1$.
Then for all but countably many~$A\in\adA$ there are $D\subseteq\N$ and a
function $H:D\to A$ such that $\phi(x)=H\inv[x]^*$ for all~$x\subseteq A$.
\end{theorem}

This theorem will now be applied to establish the following proposition.

\begin{proposition}[$\OCA$]\label{lemma.locally.complete}
For every $f\in\functions$ there is a $\beta<\omega_1$ such
that for every $\alpha\ge\beta$ there is a lifting
$\Phi_{f,\alpha}$ of~$\phi$ with $\Phi_{f,\alpha}(O)=\Phi_{f,\alpha}[O]$
whenever $O\subseteq B_{f,A_\alpha}$ and such that
$\Phi_{f,\alpha}(n,s)\cap\Phi_{f,\alpha}(m,t)=\emptyset$
whenever $\{n\}\times[s]$ and $\{m\}\times[t]$ are disjoint.
\end{proposition}

\begin{proof}
We fix $f\in\functions$ and show how to find $\beta$ and $\Phi_{f,\alpha}$
for each~$\alpha\ge\beta$.
We transfer the almost disjoint family~$\adA$ to~$B_f$ by setting
$C_\alpha=B_{f,A_\alpha}$ and $\calC=\{C_\alpha:\alpha<\omega_1\}$.

It is fairly straightforward to show that $\calC$~is neat; one stretches
the bijection~$\pi$ to find an injection~$\tilde\pi$ from~$B_f$ to~$\bintree$
that maps every~$C_\alpha$ onto a branch of~$\bintree$ and
different~$C_\alpha$ to different branches.

This means that we can apply Theorem~\ref{thm.powN.into.powNfin}
to the embedding $\phi_f$ of~$\pow(B_f)$
into~$\powNfin$, defined by $\phi_f(O)=\phi(O)$.
This gives us a~$\beta$ and for every~$\alpha\ge\beta$ a subset~$D_\alpha$
of~$\N$ and a function $H_\alpha:D_\alpha\to C_\alpha$ such that
for every subset~$O$ of~$B_f$ the set~$H_\alpha\inv[O]$ is a representative
of~$\phi_f(O)$.
We can define $\Phi_{f,\alpha}$ by $\Phi_{f,\alpha}(O)=H\inv[O]$
for $O\subseteq B_f$
and by setting $\Phi_{f,\alpha}(a)=\Psi(a)$ for the other elements of~$\M$.
\end{proof}

For each $f$ we denote the minimum possible~$\beta$ by~$\alpha_f$.

\begin{lemma}\label{lemma.alpha.monotone}
If $f\almle g$ then $\alpha_f\le\alpha_g$.
\end{lemma}

\begin{proof}
Let $\alpha\ge\alpha_g$ and consider the lifting~$\Phi_{g,\alpha}$.
We define a lifting $\Phi_{f,\alpha}$ in a fairly obvious way:
first fix~$m$ such that $f(n)\le g(n)$ for $n\ge m$
and, if need be, redefine, for the duration of this proof, the
values~$\Psi(n,s)$ for $n<m$ and~$s\in\sbintree$ so as to
get $\Psi(n,s)\cap\Phi_{g,\alpha}(l,t)=\emptyset$ whenever this is needed.

Then, given $O\subseteq B_{f,A_\alpha}$ put
$O_1=\bigl\{\orpr{n}{s}\in O:n\ge m\bigr\}$ and put
$$
U_O=\bigl\{\orpr{n}{t}\in B_{g,A_\alpha}:
            (\exists\orpr{n}{s}\in O_1)(s\subseteq t)\bigr\}.
$$
We define, for $O\subseteq B_{f,A_\alpha}$,
$$
\Phi_{f,\alpha}(O)=\Phi_{g,\alpha}[U_O]\cup
      \Psi\bigl[\bigl\{\orpr{n}{s}\in O:n<m\bigr\}\bigr].
$$
Note that we implicitly defined
$\Phi_{f,\alpha}(n,s)=
 \Phi_{g,\alpha}\bigl[\bigl\{\orpr{n}{t}\in B_g:s\subseteq t\bigr\}\bigr]$
whenever~$n\ge m$
and $\Phi_{f,\alpha}(n,s)=\Psi(n,s)$ when~$n<m$.
It follows that
$\Phi_{f,\alpha}(O_1)=\Phi_{g,\alpha}[U_O]=\Phi_{f,\alpha}[O_1]$
and hence that
$$
\Phi_{f,\alpha}(O)=\Phi_{f,\alpha}[O].
$$
We already took care of the disjointness requirement so this~$\Phi_{f,\alpha}$
witnesses that $\alpha\ge\alpha_f$ (once we use $\Psi$
to define~$\Phi_{f,\alpha}$ on the rest of~$\M$).
\end{proof}

We apply Lemma~\ref{lemma.eventually.constant} to find $\alpha_\infty$
such that $\alpha_f\le\alpha_\infty$ for all~$f$.

We fix $\alpha\ge\alpha_\infty$ and put $A=A_\alpha$.
For every $f\in\calF$ we simply write $\Phi_f$ for the
lifting~$\Phi_{f,\alpha}$.
For every~$f$ we extend $\Phi_f$ in a natural way to the set
$C_f=\bigl\{\orpr{n}{s}:n\in A,\card{s}\le f(n)\bigr\}$:
we demand that
$\Phi_f(s,n)=\Phi_f(n,\szer)\cup\Phi_f(n,\sone)$
whenever appropriate; this makes $\Phi_f$ exact on~$C_f$.

\begin{lemma}\label{lemma.Phi's.cohere}
If $f\almle g$ then $\Phi_g\restr C_f\almeq \Phi_f$.
\end{lemma}

\begin{proof}
Consider a potential sequence $\bigl<\orpr{n_i}{s_i}:i\in\omega\bigr>$
of points in~$C_f$ where $\Phi_g$ and $\Phi_f$ disagree.
By the disjointness condition and because the symmetric difference
of $\Phi_f(n_i,s_i)$ and $\Phi_g(n_i,s_i)$ is always finite we can assume
that $\Phi_f(n_i,s_i)$ does not meet~$\Phi_g(n_j,s_j)$ when~$i<j$.
Let $O_g=\bigl\{\orpr{n_i}{s}\in B_g:i\in\omega, s_i\subseteq s\bigr\}$
and $O_f=\bigl\{\orpr{n_i}{s}\in B_f:i\in\omega, s_i\subseteq s\bigr\}$.

Observe that $O_g$ and $O_f$ determine the same open subset of~$A\times\Cantor$,
so that $\Phi_g(O_g)\almeq\Phi_f(O_g)$.
It should be clear however that by the choice of the points~$\orpr{n_i}{s_i}$
we have $\Phi_g[O_g]\almneq\Phi_f[O_f]$, which is a contradiction.
\end{proof}

Observe that because $\bound=\aleph_2$ the family $\{C_f:f\in\functions\}$
is a $P_{\aleph_1}$-ideal on $A\times\bintree$;
we can therefore apply Theorem~\ref{cor.uniformizing} to find one
map~$\Phi$ from~$A\times\bintree$ to~$\powN$ such that
$\Phi\restr C_f\almeq\Phi_f$ for all $f\in\functions$.
This function~$\Phi$ is almost as required.

First choose $m\in\omega$ and a $\almle$-cofinal family~$\calF$ consisting
of increasing elements of~$\functions$ such that $\Phi(n,s)=\Phi_f(n,s)$
whenever $f\in\calF$, $n\ge m$ and $s\in\sbintree$.
Without loss of generality the set $\{f(m):f\in\calF\}$ is unbounded
--- make $m$ a bit larger if necessary (if no larger~$m$ can be found
the family $\calF$ is not even $\almle$-unbounded).

This immediately implies that $\Phi$ is exact on~$(A\setminus m)\times\bintree$;
we simply modify~$\Phi$ slightly on $(A\cap m)\times\bintree$ to make it exact
on the whole of~$A\times\bintree$.
For all other elements $a$ of~$\M$ we put $\Phi(a)=\Psi(a)$.

The proof that $\Phi$~is complete is much like the proof of
Lemma~\ref{lemma.alpha.monotone}.
Let $f\in\functions$ and $l\ge m$ such that $\Phi(n,s)=\Phi_f(n,s)$
whenever $n\ge l$.
Given $O\subseteq B_f$ we first note that $\Phi(O)\almeq\Phi_f(O)$, because
both $\Phi$ and $\Phi_f$ are liftings.

To complete the proof we show that also $\Phi[O]\almeq\Phi_f[O]$.
Indeed, let $O'=\bigl\{\orpr{n}{s}\in O:n\ge l\bigr\}$,
then $\Phi[O']=\Phi_f[O']$, so we are left with showing
$\Phi[O'']\almeq\Phi_f[O'']$, where $O''=O\setminus O'$.
But $O''$~is finite so that by Lemma~\ref{lemma.Phi.fin.add} we have
$\Phi[O'']\almeq\Phi(O'')$ and $\Phi_f(O'')\almeq\Phi_f[O'']$;
the equality $\Phi(O'')\almeq\Phi_f(O'')$ holds because both maps are
liftings.

\section{Embedding $\powN$ into $\powNfin$}
\label{sec.powN.into.powNfin}

In this section we
prove Theorem~\ref{thm.powN.into.powNfin}, thus completing
the argument for our main result.
We are given an an embedding $\phi$ of $\powN$ and a neat almost disjoint
family~$\adA=\{A_\alpha:\alpha<\omega_1\}$.
We have to find an~$\alpha$ such that for every $\beta\ge\alpha$
there are $D\subseteq\N$ and $H:D\to A_\alpha$ such that
$\phi(x)=H\inv[x]^*$ for all subsets~$x$ of~$A_\alpha$.

We begin by taking a lifting~$\Phi:\powN\to\powN$ of~$\phi$.
We may assume, upon replacing $\sPhin$ by
$\bigl(\{n\}\cup\sPhin\bigr)\setminus\bigcup_{i<n}\sPhin[i]$,
that the~$\sPhin$ form a partition of~$\N$.
We shall identify $\N$ with $\N\times\N$ in such a way that $\sPhin$
corresponds to the vertical line~$\{n\}\times\N$; we shall therefore write~$V_n$
for~$\sPhin$.

For an $f\in\functions$ we write
$L_f=\bigl\{\orpr{n}{m}:n\in\omega,m\le f(n)\bigr\}$.
The following lemma will be useful toward the end of the proof.

\begin{lemma}\label{lemma.Phi.almost.works}
For each $a\subseteq\N$ there is $f\in\functions$ such that
$\Phi(a)\setminus L_f=\bigcup_{n\in a}V_n\setminus L_f$.
\end{lemma}

\begin{proof}
If $n\in a$ then $V_n\almsubseteq\Phi(a)$ and if $n\notin a$
then $V_n\cap\Phi(a)\almeq\emptyset$; now let $f$ code witnesses:
if $n\in a$ then $V_n\setminus L_f\subseteq\Phi(a)$ and if
$n\notin a$ then $V_n\setminus L_f\cap\Phi(a)=\emptyset$.
\end{proof}

We enumerate $\adA$ as $\{A_\alpha:\alpha\in\omega_1\}$.

For $f\in\functions$ consider $\Phi_f:\powN\to\pow(L_f)$, defined by
$\Phi_f(a)=\Phi(a)\cap L_f$ and observe that $\Phi_f$ induces a
homomorphism from~$\powNfin$ to $\powNfin$.

As in~\cite{DowHart99} $\OCA$~may now be applied to give us
an~$\alpha_f<\omega_1$ such that $\Phi_f$~is simple on~$A_\alpha$
whenever~$\alpha\ge\alpha_f$, where `simple' means that there
are $D\subseteq L_f$ and a finite-to-one function $h:D\to A_\alpha$ such that
$\Phi_f(a)\almeq h\inv[a]$ for all subsets~$a$ of~$A_\alpha$.
As in the previous section we choose~$\alpha_f$ as small as possible
and we use the following lemma to fix~$\alpha_\infty$ such that
$\alpha_f\le\alpha_\infty$ for all~$f$.

\begin{lemma}
If $f\almle g$ then $\alpha_f\le\alpha_g$.
\end{lemma}

\begin{proof}
Take $\alpha\ge\alpha_g$ and fix $D\subseteq L_g$ and $h:D\to A_\alpha$
such that $\Phi_g(a)\almeq h\inv[a]$ for all $a\subseteq A_\alpha$.
Now simply let $D_1=D\cap L_f$ and $h_1=h\restr D_1$; clearly
$\Phi_f(a)\almeq\Phi_g\cap L_f\almeq h\inv[a]\cap L_f=h_1\inv[a]$
for all~$a\subseteq A_\alpha$.
We see that $\alpha\ge\alpha_f$.
\end{proof}

For the rest of the proof we fix an~$\alpha\ge\alpha_\infty$ and show that
$A=A_\alpha$~is as required.
For each $f\in\functions$ we take $D_f$ and $h_f:D_f\to A$ as above.
We intend to find $D$ and~$H$ by an application of
Theorem~\ref{thm.uniformizing}.

\begin{lemma}
If $f\le g$ then $D_f\almeq D_g\cap L_f$ and $h_g\restr D_f\almeq h_f$.
\end{lemma}

\begin{proof}
The first equality is clear: by construction $D_g\almeq\Phi_g(A)$
and $D_f\almeq \Phi_f(A)$, so that
$D_g\cap L_f\almeq \Phi_g(A)\cap L_f\almeq\Phi_f(A)\almeq D_f$.

To prove the second equality let $x$ be an infinite subset of $D_g\cap D_f$
such that $h_f(i)\neq h_g(i)$ for all~$i\in x$; because $h_f$ and~$h_g$
are finite-to-one we can assume that $h_f[x]\cap h_g[x]=\emptyset$.
But then we would have a contradiction because on the one hand
$\Phi_g\bigl(h_f[x]\bigr)\cap\Phi_g\bigl(h_g[x]\bigr)
\almeq h_g\inv\bigl[h_f[x]\bigr]\cap h_g\inv\bigl[h_g[x]\bigr]=\emptyset$
while on the other hand
$x\subseteq h_f\inv\bigl[h_f[x]\bigr]\cap h_g\inv\bigl[h_g[x]\bigr]
  \almsubseteq\Phi_g\bigl(h_f[x]\bigr)\cap\Phi_g\bigl(h_g[x]\bigr)$.
\end{proof}

We apply Theorem~\ref{thm.uniformizing} to the family $\{F_f:f\in\functions\}$
of functions defined by $F_f(p)=\bigl<\chi_{D_f}(p),h_f(p)\bigr>$
to find a function $F:\omega\times\omega\to 2\times\omega$ that uniformizes
this family.
We set $D=\{p:F_1(p)=1\}$ and $H=F_2\restr D$.

For the rest of the proof the letters $n$ and $m$ will refer to elements of~$A$.

\begin{lemma}
There are only finitely many $m$ in~$A$ for which the set
$\{n\in A:H\inv(n)\cap V_m\neq\emptyset\}$ is infinite.
\end{lemma}

\begin{proof}
Let $b$ be the set of~$m$ in~$A$ for which
$I_m=\{n\in A:H\inv(n)\cap V_m\neq\emptyset\}$ is infinite.
Thin out~$b$ to get $I_m\setminus b$ infinite for all~$m$ in~$b$.
Choose $f\in\functions$ as per Lemma~\ref{lemma.Phi.almost.works}
for~$A$ and~$b$, so that $\Phi(b)\setminus L_f=(b\times\omega)\setminus L_f$
and
$\Phi(A\setminus b)\setminus L_f=
 \bigl((A\setminus b)\times\omega\bigr)\setminus L_f$.

Because the sets $H\inv(n)$ are pairwise disjoint we can a one-to-one choice
function $m\mapsto n_m$ for the family $\{I_m\setminus b:m\in b\}$ such that
$H\inv(n_m)\cap V_m\setminus L_f\neq\emptyset$ for all~$m$.
We choose a function $g>f$ that captures these intersections:
$H\inv(n_m)\cap V_n\cap(L_g\setminus L_f)\neq\emptyset$ for all~$m$.
It follows that
$H\inv[b]\cap(L_g\setminus L_f)$ meets $(A\setminus b)\times\omega$ in an
infinite set.
However, by the choice of~$f$ and the properties of~$h_g$ the set
$H\inv[b]\cap (L_g\setminus L_f)$ is almost equal to
$\bigcup_{m\in b}\{m\}\times\bigl(f(m),g(m)\bigr]$
which is disjoint from $(A\setminus b)\times\omega$.
\end{proof}

\begin{lemma}
For every $n\in A$ the set $\{m\in A:H\inv(n)\cap V_m\neq\emptyset\}$ is
finite.
\end{lemma}

\begin{proof}
Fix $f\in\functions$ such that for all~$n$ and~$m$ in~$A$:
\emph{if} $n\neq m$ and $H\inv(n)\cap V_m\neq\emptyset$
\emph{then} $H\inv(n)\cap V_m\cap L_f\neq\emptyset$.
Now note that $H\inv(n)\cap L_f\almeq h_f\inv(n)$, so that $H\inv(n)\cap L_f$~is
finite.
\end{proof}

\begin{lemma}
There are only finitely many pairs $\orpr{n}{m}$ for which $H\inv(n)\cap V_m$ is
infinite.
\end{lemma}

\begin{proof}
Assume we have $\bigl\{\orpr{n_i}{m_i}:i\in\omega\bigr\}$ with
$H\inv(n_i)\cap V_{m_i}$ infinite for all~$i$ and $n_i,m_i<n_j,m_j$ whenever
$i<j$.
Choose $f\in\functions$ as per Lemma~\ref{lemma.Phi.almost.works} for the sets
$b=\{n_i:i\in\omega\}$ and $c=\{m_i:i\in\omega\}$ and choose $g>f$ such that
$H\inv(n_i)\cap V_{m_i}\cap(L_g\setminus L_f)\neq\emptyset$ for all~$i$.

We obtain a contradiction as before: $b$ and~$c$ are disjoint, hence $\Phi(b)$
and $\Phi(c)$ are almost disjoint.
On the other hand $h_g\inv[b]\setminus L_f\almsubseteq\Phi(b)$
and $\bigcup_{m\in c}\{m\}\times\bigl(f(m),g(m)\bigr]\subseteq\Phi(c)$;
the intersection of the smaller sets is infinite.
\end{proof}

Putting these lemmas together we see that there are $M$ and $N$ in~$\omega$
such that $N\ge M$ and
\begin{enumerate}
\item if $n,m\ge M$ and $n\neq m$ then $H\inv(n)\cap V_m$ is finite;
\item if $m\ge M$ then $V_m$ meets only finitely many~$H\inv(n)$; and
\item if $n<M$ and $m\ge N$ then $H\inv(n)\cap V_m=\emptyset$.
\end{enumerate}
(Note that $M$ should be chosen first, to ensure 1 and 2.)

By 1 and 2 we can fix $h\in\functions$ such that
$H\inv(n)\cap V_m\subseteq L_h$
whenever $n,m\ge M$ and $n\neq m$; an application of~3 then tells us that
$H\inv(n)\setminus L_h=V_n\setminus L_h$ for $n\ge N$.
We also redefine $H$ on the set $N\times\omega$ to get $H\inv(n)=V_n$ for $n<N$.

Now let $b\subseteq A$ and fix $f>h$ as per Lemma~\ref{lemma.Phi.almost.works}.
By the choice of~$f$ and~$h$ we have
$$
\Phi(b)\setminus L_f=(b\times\omega)\setminus L_f=H\inv[b]\setminus L_f.
$$
The redefined~$H$ still satisfies $H\restr L_f\almeq h_f$ so that
$$
\Phi(b)\cap L_f\almeq h_f\inv[b]\almeq H\inv[b]\cap L_f.
$$
This shows that $H$~is as required.

\section{Comments and questions}
\label{sec.comments}

In this section we comment on two earlier results about the embeddability
of~$\M$ into~$\powNfin$ and raise some questions.

\subsection{Martin's Axiom does not imply embeddability}

A consequence of our result is that Martin's Axiom does not imply that
$\M$~can be embedded into~$\powNfin$; this is so because the
conjunction $\OCA+\MA$ is consistent --- it follows from the Proper
Forcing Axiom and it can be proved consistent in an $\omega_2$-length
iterated forcing construction.

In \cite{FrankiewiczGutek81} Frankiewicz and Gutek assert that $\MA$
implies there is a measure-preserving embedding~$\phi$ of~$\M$
into~$\powNfin$ --- measure preserving in the sense that for every
element~$a$ one has $\mu(a)=d\bigl(\phi(a)\bigr)$.
Here $d$~is the \emph{asymptotic density}, defined by
$$
d(X)=\lim_{n\to\infty}\frac{\bigcard{X\cap\{1,2,\ldots,n\}}}n,
$$
for those subsets~$X$ of~$\N$ for which the limit exists.
Of course for this to make sense we must consider the standard
incarnation of~$\M$ as the quotient of the Borel algebra of the unit interval
by the ideal of measure-zero sets.

The reader is likely to be interested where the argument has a gap.
It seems that the principal error in their proof is in the following
lemma, which is the key step in the construction of the embedding.

\begin{lemma}[$\MA$]\label{lemma.f-g.false}
If\/ $\LL$ is a subalgebra of~$\M$ of size less than~$\cont$ and if
$a\in\M\setminus\LL$ then every measure-preserving embedding~$\phi$
of\/~$\LL$ into~$\powNfin$ can be extended to a measure-preserving
embedding~$\psi$ of the algebra generated by\/ $\LL\cup\{a\}$
into~$\powNfin$.
\end{lemma}

It is relatively easy to see that this lemma is true, in $\ZFC$, for
countable~$\LL$ --- indeed, a value~$\psi(a)$ is readily constructed
by recursion.
This lemma is false for $\LL$ of size~$\aleph_1$, as can be seen from
the following example.
We work on the interval~$(0,1]$.
We split $(0,1]$ into intervals~$\{I_n:n\in\N\}$, where
$I_{2n-1}=(2^{-(n+1)},2^{-n}]$ and $I_{2n}=I_{2n-1}+\frac12$.
Next let $\adA=\{A_\alpha:\alpha<\omega_1\}$ be an $\aleph_1$-sized
almost disjoint family, where $A_\alpha$ consists of even/odd numbers
whenever $\alpha$~is even/odd.
For $\alpha<\omega_1$ we put $a_\alpha=\bigcup_{n\in A_\alpha}I_n$;
note that $\{a_\alpha:\alpha<\omega_1\}$ is an independent family in~$\M$
in that no element is in the algebra generated by the other elements of the
family --- we let $\LL$ be the subalgebra of~$\M$ generated by this family.
Observe that $a=(0,\frac12]$ does not belong to~$\LL$.

Let $\calL=\{L_\alpha:\alpha<\omega_1\}$ be a Luzin-type almost disjoint
family, which means that for no two disjoint uncountable subsets~$S$ and~$T$
of~$\omega_1$ one can find a subset~$X$ of~$\N$ such that
$L_\alpha\almsubseteq X$ for $\alpha\in S$
and $L_\alpha\cap X\almeq\emptyset$ for $\alpha\in T$ ---
see~\cite{Luzin47} or~\cite[Theorem~4.1]{vanDouwen84} .
We assume the family~$\calL$ lives on the set~$Z=\{n!:n\in\N\}$ ---
note that $d(Z)=0$.

We can apply the countable version of Lemma~\ref{lemma.f-g.false} to find
a measure-preserving embedding from the algebra~$\LL$
into~$\pow(\N\setminus Z)/\fin$.
Now augment this embedding to obtain an embedding~$\phi$ of~$\LL$
into~$\powNfin$ such that $\phi(a_\alpha)\cap Z=L_\alpha$ for all~$\alpha$.
There is no way to extend~$\phi$ to an embedding of~$\LL\cup\{a\}$
into~$\powNfin$, measure-preserving or otherwise:
because $a_\alpha<a$ if $\alpha$~is odd and $a_\alpha\wedge a=0$ if
$\alpha$~is even we should have $L_\alpha\almsubseteq \psi(a)$ for
odd~$\alpha$ and $L_\alpha\cap\psi(a)\almeq\emptyset$ for even~$\alpha$,
which is impossible by the choice of~$\calL$.

\subsection{The oracle-c.c.\ method}

In \cite{Frankiewicz85} Frankiewicz sketched a proof of the consistency
of the nonembeddability of~$\M$, using Shelah's oracle-c.c.\ method.
Certainly this sketch is quite incomplete but it also appears to be
fundamentally  erroneous.
Before we can point out where Frankiewicz went wrong we give a rough
outline of the oracle-c.c.\ method.

Shelah's oracle-c.c.\ method --- see
\cite{Burke93}, \cite[Chapter~IV]{Shelah82} and \cite[Chapter~IV]{Shelah98}
--- provides us with a way of preventing that various undesirable
objects will appear once an iterated forcing construction is underway.
It would take us too far afield to explain this method in full, suffice it
to say that one builds c.c.c.\ partial orders of size~$\aleph_1$ whose
antichains are kept under tight control by a $\diamond$-sequence
(an oracle).

A consistency proof, along these lines, for the nonembeddability of~$\M$
would run as follows:
given an embedding $\phi:\M\to\powNfin$, construct a partial order~$\Po$
and a $\Po$-name~$\dot X$ for a Borel set such that there is \emph{no}
$\Po*\Cohen$-name~$\dot A$ for a subset of~$\N$ that satisfies the
following two requirements for every Borel set~$Y$ from the ground
model:
$$
{}\forces
\mbox{``\,$Y\subseteq\dot X\Rightarrow\phi(Y)\almsubseteq\dot A$\,''}
\text{ and}
{}\forces
\mbox{``\,$\dot X\subseteq Y\Rightarrow\dot A\almsubseteq\phi(Y)$\,''}
$$
The extra factor $\Cohen$ is a necessary technical device
(implicit in any finite-support iteration)
that will enable one to improve a given $\diamond$-sequence into an oracle
with the property that if the rest of the iteration is kept under its
control one will never encounter an undesirable name~$\dot A$ as above.
This, together with a reflection argument involving $\diamond$
on~$\omega_2$, will ensure that ultimately no embedding of~$\M$
into~$\powNfin$ will remain.

Frankiewicz' sketch is based on Shelah's proof, from~\cite{Shelah:185},
that it is consistent that the natural quotient homomorphism~$q$
from~$\Bor(\C)$ onto~$\M$ does not split, i.e., there is no homomorphism
$h:\M\to\Bor(\C)$ such that $q\circ h=\operatorname{Id}_\M$
(such a homomorphism is also called a lifting of~$\M$ into~$\Bor(\C)$).

The sketch consists basically of two parts.
\begin{enumerate}
\item A copy of Shelah's construction of the partial order
      from~\cite{Shelah:185} with some (questionable) modifications, and
\item the unsupported assertion that
      ``From now on the proof goes exactly as in Shelah~\cite{Shelah:185}''.
\end{enumerate}
The problem with the second part is that Shelah's argument was written up in
such a way that it would also apply to the category algebra; this is the
quotient algebra~$\Kat=\Bor(\C)/\Meag$, where $\Meag$~is the ideal of meager
sets.
As none of the modifications is $\M$-specific Frankiewicz' sketch would also
lead to a proof that it is consistent for there to be no embedding of the
category algebra into~$\powNfin$.
This, however, is impossible: the category algebra is known to be
embeddable into~$\powNfin$.
This is most readily seen by noting its Stone space is separable and hence
a continuous image of~$\betaNminN$.

One very problematic modification occurs almost at the beginning:
we are given a countable partial order~$\Po$ and a
$\Po*\Cohen$-name~$\dot A$ for a subset of~$\N$.
We are then promised an extension~$\Po'$ of~$\Po$ and an infinite
subset~$B$ of~$\N$ such that some condition in~$\Po'*\Cohen$ will force
either $B\subseteq\dot A$ or $B\cap\dot A=\emptyset$.
This is patently impossible if $\dot A$ happens to be the name of the
generic subset that is added by~$\Cohen$: it is well-known that this
set and its complement meet every infinite subset of~$\N$ from the ground
model.

\subsection{Some questions}

It would be of interest to know whether there are other ways of
proving the nonembeddability of the measure algebra consistent.
To be specific we ask.

\begin{question}
Can the consistency of the nonembeddability of~$\M$ be established by
the oracle-c.c.\ method?
\end{question}

Because the oracle-c.c.\ method relies on having $\diamond$ in every
intermediate model it seems to produce models with $\cont\le\aleph_2$ only;
the known models for~$\OCA$ satisfy $\cont=\aleph_2$ as well.
Therefore the following question `really' asks for a new method.

\begin{question}
Is the nonembeddability of~$\M$ consistent with larger values for~$\cont$?
\end{question}

Many statements that were proved consistent by the oracle-c.c.\ method were
later shown to follow from $\OCA$ (or $\OCA+\MA$), see, e.g.,
\cite{Just92b}.
One of the major questions that remains is.

\begin{question}
Does $\OCA$ or even~$\PFA$ imply that the quotient homomorphism
from~$\Bor(\C)$ onto~$\M$ (or onto~$\Kat$) does not split?
\end{question}

Another question is to delineate what the subalgebras of $\powNfin$
(or, dually, the zero-dimensional continuous images of $\betaNminN$) look
like.
The $\ZFC$ results can, broadly, be divided into two categories:
the first contains the easy result that every separable compact
space is a continuous image of $\betaNminN$;
the second contains moderately difficult constructions of embeddings
into~$\powNfin$ of assorted algebras --- these are generally small in some sense
or another, so that the construction can be seen to terminate.
The smallness conditions all seem to imply that the algebra in question is
incomplete.
This suggests that the following question might have a positive answer.

\begin{question}
Is it consistent that every complete Boolean algebra that is embeddable
into~$\powNfin$ must be $\sigma$-centered?
\end{question}


\providecommand{\bysame}{\leavevmode\hbox to3em{\hrulefill}\thinspace}

\end{document}